\documentclass[11pt,letterpaper]{article}

\usepackage{etex}
\usepackage{amsmath}
\usepackage{amsfonts}
\usepackage{amssymb}
\usepackage{multirow}
\usepackage{multicol}
\usepackage{longtable}
\usepackage[usenames,dvipsnames]{color}
\usepackage[top=1in, bottom=2in, left=1in, right=1in]{geometry} 
\usepackage{pdflscape}
\usepackage{graphicx}
\usepackage{setspace}
\usepackage{palatino}
\usepackage{url}
\usepackage[pdftex,breaklinks]{hyperref}

\newcommand{\bit}{\begin{itemize}}
\newcommand{\eit}{\end{itemize}}
\newcommand{\ben}{\begin{enumerate}}
\newcommand{\een}{\end{enumerate}}
\newcommand{\tm}{\textrm}
\newcommand{\mm}{\mathrm}

\newcommand{\T}{\rule[2.6ex]{0pt}{0pt}}
\newcommand{\B}{\rule[-1.2ex]{0pt}{0pt}}

\hypersetup{
    unicode=false,          
    pdftoolbar=true,        
    pdfmenubar=true,        
    pdffitwindow=false,     
    pdfstartview={Fit},     
	pdfpagelayout={SinglePage},		
    pdftitle={},    							
    pdfauthor={Anshul Agarwal},		
    pdfsubject={Documentation},		
    pdfcreator={Anshul Agarwal},	
    pdfproducer={Anshul Agarwal},	
    pdfkeywords={}, 							
    pdfnewwindow=true,      
    colorlinks=true,        
    linkcolor=blue,         
    citecolor=blue,         
    filecolor=blue,         
    urlcolor=blue           
}

\author{Anshul Agarwal\thanks{Corresponding author, anshul.agarwal@united.com, 
United Airlines, Chicago, IL 60606, USA}}
\date{}
\title{Validation of Inventory models for Single-echelon Supply Chain using Discrete-event Simulation}

\begin{document}

\maketitle

\begin{abstract}

Inventory decision frameworks proposed in the literature for single-echelon supply chain systems rely on assumptions to obtain closed form expressions.  In particular, two such frameworks - one conventional and the other with a demand undershoot - determine optimal reorder point for a desired $\beta$ or Type-II service level.  In this work we assess the accuracy and applicability of these frameworks with the help of a discrete event simulation framework developed in SimPy.  While in several cases the closed form literature models under-predict the service level, and thus result in a higher reorder point and safety stock, we observe some situations where model predicted service levels match the simulated results with statistical significance.

\end{abstract}

\section{Introduction}
\label{sec:intro}

Various models have been proposed in the literature for single-echelon supply chain systems (which we discuss further) that help determine how much inventory should be stocked in order to meet the desired customer service level.  Such models estimate the parameters for an inventory policy followed by the single-echelon system.  This work primarily assesses the accuracy and applicability of these models with the help of a simulation framework developed in SimPy \cite{matloff2008introduction}.

A supply chain facilitates transportation of goods from one point to another is a network of facilities, such as suppliers, production locations, distribution centers.  This can involve procurement of material, transformation of material to intermediate and finished products, and distribution of finished products to customers.  The operation of a supply chain is driven by customer orders for finished products.  It may be desired to stock inventories at various locations in a supply chain in order to achieve smoother operation and customer satisfaction.  However, optimal inventory levels may be imperative for a cost effective operation and reduction in working capital \cite{chopra2007supply, ganeshan1999managing, hadley1963analysis, evant2003, lee1993material, zipkin2000foundations}.

Inventory decisions such as safety stock, replenishment size, average on-hand inventory, carrying costs or working capital, etc.\ that meet a certain customer service level are challenging and not intuitive.  Complexity arises due to various sources of uncertainties, such as demand volume, lead time (process uncertainties, machine downtime, transportation reliability), supplier reliability (raw material availability, issues at the supplier end), etc.  Several frameworks have been designed in the literature that estimate aforementioned inventory decisions while accounting for such sources of uncertainties.  However, they rely on assumptions to obtain analytic closed form expressions.  In this work we validate the accuracy and applicability of two such frameworks given their respective assumptions with the help of a more accurate discrete-event simulation framework.  Note that this work deals with only the single-echelon systems.  Its extension for multi-echelon networks will be covered in future articles.

The article is organized as follows.  We first describe the single-echelon inventory management system considered in this work.  Section 3 illustrates the inventory optimization models that we validate.  In Section 4 we present the SimPy discrete-event simulation model, while the validation results are documented in the next section.  Section 6 concludes the article.

\section{Problem Scope}

\subsection{Inventory System}

As mentioned above, we consider a single-echelon system in this work.  The term \emph{echelons} for a supply chain refers to the number of interacting layers of inventory stocking locations before the product is delivered to the final customer.  In a single-echelon system, there exists a single layer of locations that receive replenishment and deliver customer shipments.  With two or more layers, the supply chain becomes multi-echelon.  In a multi-echelon system, interaction between layers can exist in a variety of forms.

The following definitions are used to define and quantify inventory in a supply chain \cite{jacobs2008operations}:
\bit
	\item \textbf{Lot size or Economic Order Quantity (EOQ)} \\
	The amount requested in each replenishment order is called \emph{lot size}.  While it may be advantageous to place several orders with a tiny lot size to lower inventory holding costs, it may not be practical due to either a fixed cost for placing an order or ability of the supply chain to fulfill orders in specific lot sizes, such as one full truck load (FTL), one railcar, one vessel shipment, etc.  The lot size that minimizes total holding and order placement cost while accounting for such system restrictions is called \emph{economic order quantity (EOQ)}.	
	\item \textbf{Cycle Stock} \\
	The time gap between initiating a replenishment order and receiving the shipment is called \emph{lead time}.  Lead time can involve pre-processing time (place order on the books, other bookkeeping), order processing time (manufacturing, assembling, product wheel, scheduling, etc.), and post processing (transportation, load-unload, etc.).  The inventory maintained to serve customers orders received during the lead time is called \emph{cycle stock}.
	\item \textbf{Safety Stock} \\
	\emph{Safety stock} is the buffer inventory required to insure against the variability in customer demand (inaccurate forecast) and lead time (order processing delays, transit delays, etc.)
\eit

\subsection{Inventory Tracking}

\begin{figure}
	\centering
	\resizebox{16cm}{!}{\includegraphics[scale=1.0]{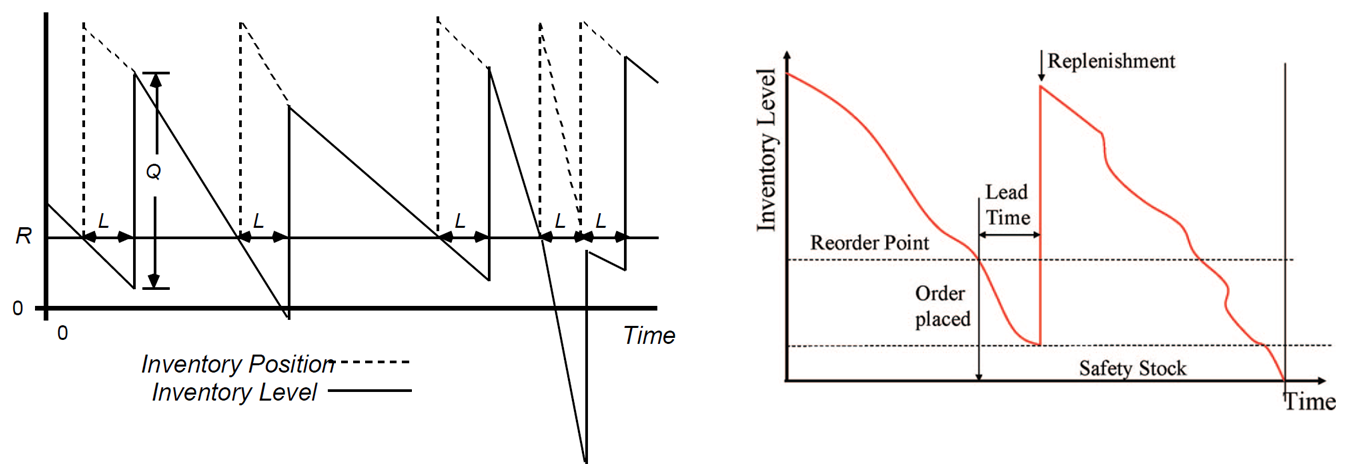}}
	\caption{Reorder point policy \cite{jacobs2008operations}}
	\label{fig:rop}
\end{figure}

Inventory policies provide a structure to manage ordering and receipt of products, timing the order placement, keeping track of what has been ordered, and deciding how much to order.  In this work we evaluate a \emph{reorder point policy}.  In a reorder point policy, as shown in Figure \ref{fig:rop}, an order is placed when the stock reaches a predetermined number of units, called the reorder point.  An order of fixed size $Q$ is placed when the inventory position (on-hand plus on-order minus backordered quantities, if any) reaches the reorder point $R$ \cite{jacobs2008operations}.  Note that the event of reaching $R$ is not ``time-triggered'' and may take place at any time depending on the demand.  The reorder point $R$ is set to cover the cycle stock plus a safety stock determined by the desired service level.  The fixed order size $Q$ is generally governed by a fixed lot size or EOQ.

\subsection{System Performance}

An inventory system's performance is measured using \emph{service level}.  The event in which the inventory is completely exhausted and customer demand for the product cannot be fulfilled is called a \emph{stockout}.  Service level is a metric to quantify the level of stockouts that can be tolerated in order to meet the desired level of customer satisfaction \cite{shivsharan2012optimizing}.  If the demand is fulfilled later, it is called a \emph{backorder}.  In such cases, the service level measures how much demand was fulfilled on time.  Inventory policies offer levers to manipulate in order to balance the trade-off between service level and the amount of inventory held.

\begin{figure}
	\centering
	\resizebox{10cm}{!}{\includegraphics[scale=1.0]{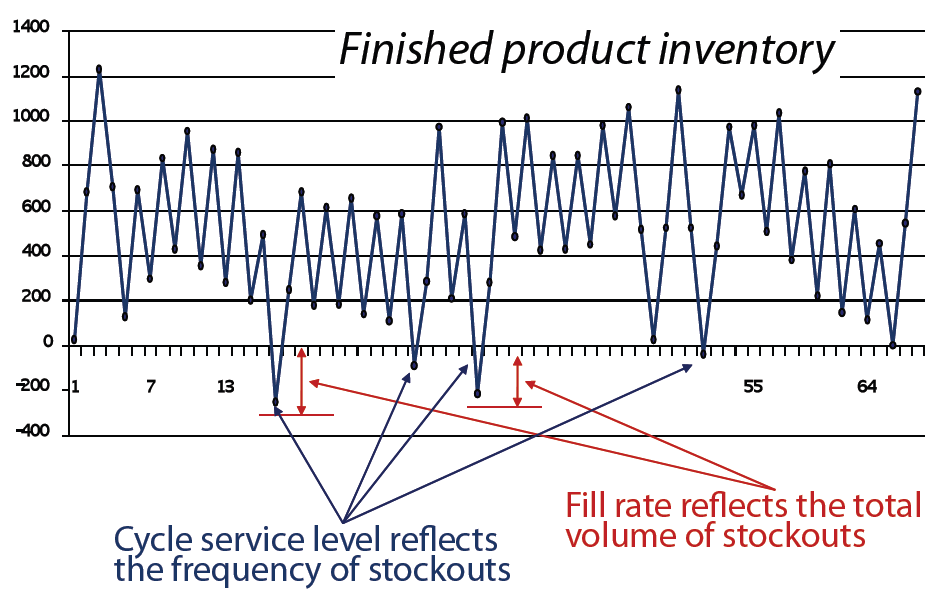}}
	\caption{Illustration of $\alpha$ and $\beta$ service level \cite{kingcrack}}
	\label{fig:service}
\end{figure}

There are two kinds of definitions of service level that are most frequently used: \emph{$\alpha$ service level (Type I)} (also called cycle service level) and \emph{$\beta$ service level (Type II)} (also known as fill rate).  The $\alpha$ service level measures the number of inventory cycles in which stockout (or backorder) did not occur.  Specifically,
\begin{equation*}
\alpha \tm{ service level} = \frac{\tm{\# of cycles with no stockouts}}{\tm{Total \# of cycles}}
\end{equation*}
It measures the probability of not stocking out during the lead time.  However, it does not measure the ``magnitude'' of stockout (if the stockout was 1 unit or 1000 units of demand).  On the other hand, the $\beta$ service level measures the fraction of the amount of demand met from the inventory (see Figure \ref{fig:service}).  Specifically,
\begin{equation*}
\beta \tm{ service level} = \frac{\tm{Total volume of orders supplied}}{\tm{Total demand ordered}}
\end{equation*}
For example, for the 28 replenishment cycles in Figure \ref{fig:service}, we observe stockouts happening in 4 cycles.  Thus, $\alpha$ service level = 24/28 = 86\%.  In contrast, the total customer demand in the 28 cycles is 16990 units, while the total stockout volume is 520 units.  Therefore, the $\beta$ service level is (16990-520)/16990 = 97\%.  In general, the $\beta$ service level tends to be higher than $\alpha$.  However, high demand and lead time variability can cause the former to be lower than the latter as the volume of stockouts go up \cite{kingcrack}.

In this work, we only focus on the models that use $\beta$ service level as a criterion to determine required inventory.

\section{Models for Validation}
\label{sec:models}

Determining optimal inventory in a reorder point system essentially involves determining the reorder quantity ($Q$ in Figure \ref{fig:rop}), the reorder point ($ROP$)($R$ in Figure \ref{fig:rop}) and the corresponding safety stock ($SS$). The value of $Q$ is obtained via minimum lot size of the system or $EOQ$ determined from ordering cost.  The models in the literature estimate only the $ROP$ and $SS$ for a given $Q$.  The choice of $Q$ impacts $ROP$ and $SS$.  The models developed in the literature are valid only when demand follows a Normal distribution.  However, no distribution assumption is required for lead time variability.

\noindent Let \\
$\mu_{D}$:  mean of the Normally distributed demand, \\
$\sigma_{D}$:  standard deviation of the Normally distributed demand, \\
$\mu_{L}$:  mean of the variable lead time, and \\
$\sigma_{L}$:  standard deviation of the variable lead time. \\
We define expected demand during lead time ($\mu$) and total deviation ($\sigma$) as
\begin{gather}
\label{eq:mean}
	\begin{array}{l}
		\mu = \mu_{D}\mu_{L}  \\
		\sigma = \sqrt{\mu_{L}\sigma_{D}^{2} + \mu_{D}^{2}\sigma_{L}^{2}}
	\end{array}
\end{gather}
The cycle stock ($CS$), $SS$, and $ROP$ are calculated as \cite{chopra2007supply,cognizantss,kingcrack}
\begin{gather}
	\begin{array}{l}
		CS = \max(Q,\mu) \\
		SS = \lambda\sigma \\
		ROP = CS + SS
	\end{array}
\end{gather}
Here $\lambda$, called the safety stock factor, varies with $\alpha$ and $\beta$ service levels \cite{chopra2007supply,cognizantss}.  We consider the following two models in this work for $\beta$ service level.

\subsection{Conventional Model}

In the conventional model, as defined in \cite{chopra2007supply,jacobs2008operations,cognizantss,kingcrack}, the $\beta$ service level is estimated as
\begin{gather}
\label{eq:fill}
\beta = 1-E/Q \Rightarrow E = (1-\beta)Q
\end{gather}
where, the \emph{expected shortage} ($E$) is the demand backordered and not fulfilled on time, and $Q$ is the reorder quantity.  The expected shortage is the excess demand during the lead time that could not be met from the cycle and safety stock. If $f(x)$ is the demand distribution, the expected shortage is calculated as
\begin{gather}
\label{eq:shortage}
E = \int_{ROP}^{\infty} (x-ROP)f(x) \mm{d}x
\end{gather}
When demand is Normally distributed, Equation (\ref{eq:shortage}) can be simplified to
\begin{gather}
\label{eq:shortnormal}
E = \sigma(\phi(\lambda) - \lambda[1-\Phi(\lambda)])
\end{gather}
Here $\sigma$ is obtained from Equation (\ref{eq:mean}), and $\phi(\lambda)$ and $\Phi(\lambda)$ are the standard normal density and the cumulative distribution functions, respectively, as below
\begin{gather}
\phi(\lambda) = (1/\sqrt{2\pi})\mm{e}^{-\lambda^{2}/2}, \quad \Phi(\lambda)=\int_{-\infty}^{\lambda}\phi(x)\mm{d}x \nonumber
\end{gather}
To approximate the integral, the following empirical expression for $\Phi(\lambda)$ can be used \cite{aludaat2008note}
\begin{gather}
\label{eq:normapprox}
\Phi(\lambda) = \frac{\mm{e}^{2y}}{1+\mm{e}^{2y}}, \quad y=0.7988\lambda(1+0.04417\lambda^{2})
\end{gather}
Equations (\ref{eq:fill}), (\ref{eq:shortnormal}), and (\ref{eq:normapprox}) can be solved together to obtain $\lambda$ for a desired $\beta$.  For the lost sales case, $(1-\beta)$ can be replaced with $(1-\beta)/\beta$ in Equation (\ref{eq:fill}).

\subsection{Undershoot Model}

In a conventional reorder point system, an order is placed when the inventory position is \emph{exactly} equal to the reorder point.  The conventional model above represents such a system.  However, this requires inventory to be monitored constantly.  In practice, inventory levels are monitored at a pre-determined frequency, such as once a day.  Consequently, a delay is observed between when $ROP$ is reached and when the order is placed.

\begin{figure}
	\centering
	\resizebox{12cm}{!}{\includegraphics[scale=1.0]{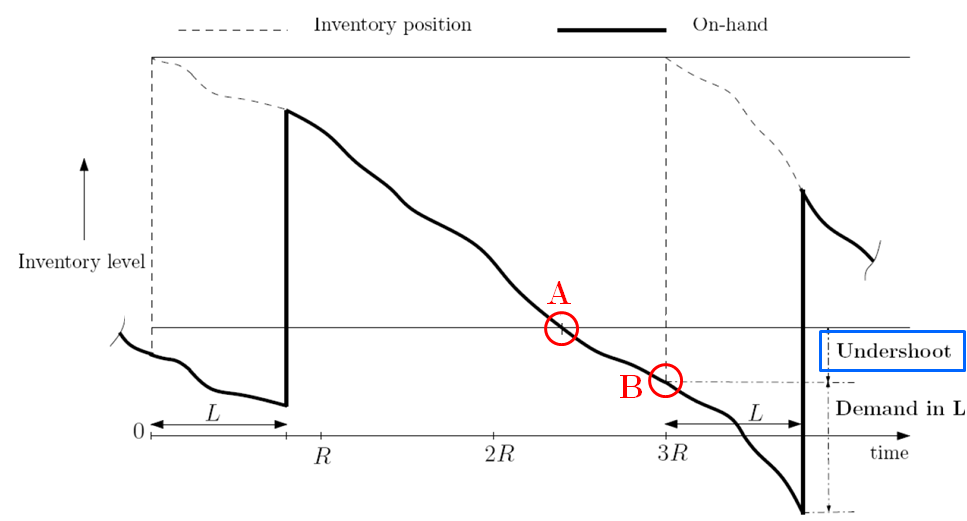}}
	\caption{The concept of reorder point undershoot \cite{Silver2008}}
	\label{fig:undershoot}
\end{figure}

Figure \ref{fig:undershoot} illustrates this concept.  Here inventory is monitored at a constant interval $R$.  The reorder point is reached at A.  While a conventional system will place an order at A, here the order gets placed at B when the system monitors inventory after $3R$.  The difference between A and B is called \emph{undershoot}.  Undershoot results in a higher $ROP$ because inventory needs to be kept not only to meet demand during the lead time but also to serve the undershoot.  As a result, expected shortage can be higher in such systems, and is not adequately captured by Equations (\ref{eq:shortage}) and (\ref{eq:shortnormal}).

Let $R$ be the inventory monitoring frequency of the system.  The $ROP$ is calculated as \cite{Silver2008}:

\begin{align}
\mu_{R} &= \mu_{D}R  &  \sigma_{R} &= \sigma_{D}\sqrt{R} \nonumber \\
\mu_{LT} &= \mu_{D}\mu_{L}  &  \sigma_{LT} &= \sqrt{\mu_{L}\sigma_{D}^{2} + \mu_{D}^{2}\sigma_{L}^{2}} \nonumber \\
ROP &= \mu_{R}+\mu_{LT}+\lambda\sqrt{\sigma_{R}^{2}+\sigma_{LT}^{2}} \label{eq:ropnew}
\end{align}
In order to calculate $\lambda$, we consider the following steps \cite{schneider1978methods,tijms1984simple}.  First we calculate the expected value of the undershoot.
\begin{align}
\label{eq:undershoot}
E_{U} &= (\mu_{R}^{2} + \sigma_{R}^{2})/2\mu_{R}
\end{align}
Next we find the expected value of shortage, which is given by
\begin{align}
E &= \frac{1}{2\mu_{R}}\int_{ROP}^{\infty} (x-ROP)^{2}f(x) \mm{d}x
\end{align}
Since demand is Normally distributed, this can be simplified to
\begin{align}
\label{eq:shortnew}
E &= \frac{\sigma_{R}^{2}+\sigma_{LT}^{2}}{2\mu_{R}}\left[ (1+\lambda^{2})\{1-\Phi(\lambda)\}-\lambda\phi(\lambda) \right]
\end{align} 
Finally, $\lambda$ is obtained by solving the following equation for a given $Q$ and $\beta$.
\begin{align}
\label{eq:findbetanew}
E &= (1-\beta)[Q + E_{U}]
\end{align}

\section{SimPy Simulation}

We developed a discrete-event simulation model in SimPy \cite{matloff2008introduction} in order to assess the accuracy and validity of the aforementioned two models.  Complete computer codes for the simulation are available here \cite{simpy}.

To simulate an inventory storage node in a single-echelon supply chain, we create a generic \emph{Process} in the SimPy simulation model called \emph{stockingFacility}.  The activation function of this process is called \emph{runOperation} that is run at a fixed event rate of one time unit.  Each \emph{runOperation} call records the random demand in that event and the shipment sent to the customer, and updates the backorder and inventory position of the node, as shown below:
\begin{verbatim}

class stockingFacility(Process):

    def runOperation(self):
        while True:
            yield hold, self, 1.0
            demand = float(np.random.normal(self.meanDemand, \
            				self.demandStdDev, 1))
            self.totalDemand += demand
            shipment = min(demand + self.totalBackOrder, \
            					self.on_hand_inventory)
            self.on_hand_inventory -= shipment
            self.inventory_position -= shipment
            backorder = demand - shipment
            self.totalBackOrder += backorder
            self.totalLateSales += max(0.0, backorder)
            if self.inventory_position <= 1.01*self.ROP:
                order = newOrder(self.ROQ)
                activate(order, order.ship(self))
                self.inventory_position += self.ROQ

\end{verbatim}
Note that when the inventory position goes below the reorder point, a new \emph{order} object is created which represents the new order placed by \emph{stockingFacility} object.  The order object is another SimPy \emph{Process} as shown below.  It's activation function \emph{ship} ensures that the replenishment is being held for a random lead time before being delivered to the \emph{stockingFacility} object.
\begin{verbatim}

class newOrder(Process):

    def ship(self, stock):
        leadTime = int(np.random.uniform(stock.minLeadTime, \
        							stock.maxLeadTime, 1))
        yield hold, self, leadTime
        stock.on_hand_inventory += self.orderQty

\end{verbatim}
The aforementioned \emph{Process} for \emph{stockingFacility} models the scenario when unfulfilled demand is backordered and is shipped at a later stage once inventory is available.  For the case when unfulfilled demand is considered lost sales, we use the following definition of the \emph{stockingFacility} \emph{Process}.  The \emph{newOrder} \emph{Process} remains same for both cases.
\begin{verbatim}

class stockingFacility(Process):

    def runOperation(self):
        while True:
            yield hold, self, 1.0
            demand = float(np.random.normal(self.meanDemand, \
            				self.demandStdDev, 1))
            self.totalDemand += demand
            shipment = min(demand, self.on_hand_inventory)
            self.totalShipped += shipment
            self.on_hand_inventory -= shipment
            self.inventory_position -= shipment
            if self.inventory_position <= 1.01*self.ROP:
                order = newOrder(self.ROQ)
                activate(order, order.ship(self))
                self.inventory_position += self.ROQ

\end{verbatim}
The simulation is initiated and activated through the following module.  We run the simulation for 1 year for all test cases.  Once the simulation is complete, the \emph{stockingFacility} object is returned with all its statistics.
\begin{verbatim}

def simulateNetwork(initialInv, ROP, ROQ, meanDemand,\ 
			demandStdDev, minLeadTime, maxLeadTime):
    initialize()
    s = stockingFacility(initialInv, ROP, ROQ, meanDemand,\
	    	demandStdDev, minLeadTime, maxLeadTime)
    activate(s, s.runOperation())
    simulate(until=365)
    s.serviceLevel = 1-s.totalLateSales/s.totalDemand #for backorder
    s.serviceLevel = s.totalShipped/s.totalDemand #for lost sales
    return s 

\end{verbatim}
The SimPy simulation framework represents the Undershoot modeling framework more closely compared to the Conventional model.  Here inventory is monitored per every one time event instead of continuous.  Therefore, by construct, we expect the simulation results to more closely match the Undershoot model vs.\ the Conventional model.

\section{Model Validation}

We validate the Conventional and Undershoot models (Section \ref{sec:models}) using the above simulation frameworks.  The model parameters are listed in Table \ref{tbl:params}.  We perform validation for three different values of demand standard deviation and four different values of the reorder quantity.

\begin{table}[b]
\begin{center}
\caption{Model parameters for the validation study}
\begin{tabular}{|l|l|}
\hline
Lead time \T & $\tilde{}Uniform(7,13)$ \\
$\lambda$ & 0 (no safety stock) \\
$ROP$ & 5000 \\ 
Mean demand & 500 \\
Demand std. deviation & [200, 400, 600] \\
$ROQ$ \B & [1000, 2000, 4000, 6000] \\
\hline
\end{tabular}
\label{tbl:params}
\end{center}
\end{table}

We run 100 replications of the simulation and collect a distribution of $\beta$ service level obtained from all replications.  The model is considered to match simulation successfully if the $\beta$ service level predicted by the model falls within two standard deviations of the $\beta$ service level distribution from the simulation.  We run simulation for both backorder and lost sales scenarios.

\begin{table}
\begin{center}
\caption{Validation for the Conventional Model - Backorder scenario}
\begin{tabular}{|c|c|c|c|c|c|}
\hline
\multirow{2}{*}{$\sigma_{D}$} & \multirow{2}{*}{$Q$} & \T Model \B & \multicolumn{2}{c|}{\T SimPy $\beta$ \B} & \T Model $\beta$ \\ \cline{4-5}
 & & $\beta$ & \T Avg. \B & Std. Dev. & matches SimPy \\
\hline
\multirow{4}{*}{200} & \T 1000 & 57.2\% & 89.3\% & 2.7\% & No \\
 & 2000 & 78.6\% & 89.7\% & 2.9\% & No \\
 & 4000 & 89.3\% & 94.8\% & 1.6\% & No \\
 & \B 6000 & 92.9\% & 95.9\% & 2.1\% & \textbf{Yes} \\
\hline
\multirow{4}{*}{400} & \T 1000 & 38.8\% & 79.5\% & 5.5\% & No \\
 & 2000 & 69.4\% & 80.3\% & 5.8\% & \textbf{Yes} \\
 & 4000 & 84.7\% & 90.1\% & 3.1\% & \textbf{Yes} \\
 & \B 6000 & 89.8\% & 89.7\% & 5.2\% & \textbf{Yes} \\
\hline
\multirow{4}{*}{600} & \T 1000 & 16.8\% & 68.1\% & 7.1\% & No \\
 & 2000 & 58.4\% & 70.3\% & 8.2\% & \textbf{Yes} \\
 & 4000 & 79.2\% & 83.7\% & 5.2\% & \textbf{Yes} \\
 & \B 6000 & 86.1\% & 81.9\% & 7.9\% & \textbf{Yes} \\
\hline
\end{tabular}
\label{tbl:convbackorder}
\end{center}
\end{table}

\begin{table}
\begin{center}
\caption{Validation for the Conventional Model - Lost sales scenario}
\begin{tabular}{|c|c|c|c|c|c|}
\hline
\multirow{2}{*}{$\sigma_{D}$} & \multirow{2}{*}{$Q$} & \T Model \B & \multicolumn{2}{c|}{\T SimPy $\beta$ \B} & \T Model $\beta$ \\ \cline{4-5}
 & & $\beta$ & \T Avg. \B & Std. Dev. & matches SimPy \\
\hline
\multirow{4}{*}{200} & \T 1000 & 70.1\% & 93.9\% & 1.2\% & No \\
 & 2000 & 82.4\% & 93.9\% & 1.2\% & No \\
 & 4000 & 90.3\% & 95.9\% & 1.2\% & No \\
 & \B 6000 & 93.3\% & 96.5\% & 1.2\% & \textbf{Yes} \\
\hline
\multirow{4}{*}{400} & \T 1000 & 62.1\% & 89.3\% & 1.9\% & No \\
 & 2000 & 76.6\% & 90.0\% & 1.8\% & No \\
 & 4000 & 86.7\% & 92.8\% & 1.8\% & No \\
 & \B 6000 & 90.7\% & 93.8\% & 2.0\% & \textbf{Yes} \\
\hline
\multirow{4}{*}{600} & \T 1000 & 54.6\% & 84.3\% & 2.7\% & No \\
 & 2000 & 70.6\% & 86.1\% & 2.5\% & No \\
 & 4000 & 82.8\% & 89.7\% & 2.4\% & No \\
 & \B 6000 & 87.8\% & 90.1\% & 2.7\% & \textbf{Yes} \\
\hline
\end{tabular}
\label{tbl:convlostsales}
\end{center}
\end{table}

Table \ref{tbl:convbackorder} shows the validation results for the conventional for the backorder scenario, while Table \ref{tbl:convlostsales} shows that for the lost sales scenario.  Here we list the average $\beta$ and its standard deviation from the SimPy simulation over 100 replications.  We also list if the model results match with simulation based on the two standard deviation rule.  Although the SimPy simulation models an inventory system with an undershoot, we use its results for comparison because the conventional reorder point system with constant inventory monitoring cannot be practically simulated.

From both tables, we observe the following:
\bit
	\item For the same set of inputs, the lost sales model exhibits a higher service level compared to the corresponding backorder model.  Moreover, the simulation variance is substantially lower for lost sales compared to backorder.
	\item For both backorder and lost sales, the model results match with simulation (within statistical limits) when $ROQ$ is higher than $ROP$.  However, for such cases, the lost sales model $\beta$ is always at least one standard deviation lower than the average simulation $\beta$.  The backorder model does not demonstrate this cadence.
	\item In the case of the backorder model, when coefficient of variation is more than 0.8, we do observe model results match with simulation when $ROQ$ is lower than $ROP$.  However, in these cases, the model $\beta$ is generally two standard deviations lower than the average simulated $\beta$.
\eit

Similarly, we present the validation results for the Undershoot model for the backorder scenario (Table \ref{tbl:ushootbackorder}) and the lost sales scenario (Table \ref{tbl:ushootlostsales}).  Here the SimPy simulation accurately models the inventory system under study.  Consequently, we expect the Undershoot model to match simulation in more cases compared to the Conventional.  However, we do not observe that in the results.  A summary of the observations is as below:
\bit
	\item Irrespective of the coefficient of variation, the model matches simulation only when $ROQ$ is more than $ROP$.
	\item The model is substantially inaccurate for cases then $ROQ$ is less than $ROP$.  Particularly for the backorder case, the model predicts absurd negative service level.
	\item When $ROQ > ROP$, for both backorder and lost sales cases, the gap between model predicted $\beta$ and average simulated $\beta$ increases as the coefficient of variation goes up.
\eit

\begin{table}
\begin{center}
\caption{Validation for the Undershoot Model - Backorder scenario}
\begin{tabular}{|c|c|c|c|c|c|}
\hline
\multirow{2}{*}{$\sigma_{D}$} & \multirow{2}{*}{$Q$} & \T Model \B & \multicolumn{2}{c|}{\T SimPy $\beta$ \B} & \T Model $\beta$ \\ \cline{4-5}
 & & $\beta$ & \T Avg. \B & Std. Dev. & matches SimPy \\
\hline
\multirow{4}{*}{200} & \T 1000 & 9.3\% & 89.3\% & 2.7\% & No \\
 & 2000 & 48.9\% & 89.7\% & 2.9\% & No \\
 & 4000 & 72.7\% & 94.8\% & 1.6\% & No \\
 & \B 6000 & 95.7\% & 95.9\% & 2.1\% & \textbf{Yes} \\
\hline
\multirow{4}{*}{400} & \T 1000 & -43.5\% & 79.5\% & 5.5\% & No \\
 & 2000 & 16.1\% & 80.3\% & 5.8\% & No \\
 & 4000 & 54.1\% & 90.1\% & 3.1\% & No \\
 & \B 6000 & 88.5\% & 89.7\% & 5.2\% & \textbf{Yes} \\
\hline
\multirow{4}{*}{600} & \T 1000 & -108\% & 68.1\% & 7.1\% & No \\
 & 2000 & -28.5\% & 70.3\% & 8.2\% & No \\
 & 4000 & 27.2\% & 83.7\% & 5.2\% & No \\
 & \B 6000 & 75.7\% & 81.9\% & 7.9\% & \textbf{Yes} \\
\hline
\end{tabular}
\label{tbl:ushootbackorder}
\end{center}
\end{table}

\begin{table}
\begin{center}
\caption{Validation for the Undershoot Model - Lost sales scenario}
\begin{tabular}{|c|c|c|c|c|c|}
\hline
\multirow{2}{*}{$\sigma_{D}$} & \multirow{2}{*}{$Q$} & \T Model \B & \multicolumn{2}{c|}{\T SimPy $\beta$ \B} & \T Model $\beta$ \\ \cline{4-5}
 & & $\beta$ & \T Avg. \B & Std. Dev. & matches SimPy \\
\hline
\multirow{4}{*}{200} & \T 1000 & 52.4\% & 93.9\% & 1.2\% & No \\
 & 2000 & 66.2\% & 93.9\% & 1.2\% & No \\
 & 4000 & 78.6\% & 95.9\% & 1.2\% & No \\
 & \B 6000 & 95.9\% & 96.5\% & 1.2\% & \textbf{Yes} \\
\hline
\multirow{4}{*}{400} & \T 1000 & 41.1\% & 89.3\% & 1.9\% & No \\
 & 2000 & 54.4\% & 90.0\% & 1.8\% & No \\
 & 4000 & 68.5\% & 92.8\% & 1.8\% & No \\
 & \B 6000 & 90.0\% & 93.8\% & 2.0\% & \textbf{Yes} \\
\hline
\multirow{4}{*}{600} & \T 1000 & 32.4\% & 84.3\% & 2.7\% & No \\
 & 2000 & 43.8\% & 86.1\% & 2.5\% & No \\
 & 4000 & 57.9\% & 89.7\% & 2.4\% & No \\
 & \B 6000 & 84.8\% & 90.1\% & 2.7\% & \textbf{Yes} \\
\hline
\end{tabular}
\label{tbl:ushootlostsales}
\end{center}
\end{table}

\section{Conclusions}

We presented two models from the literature that predict the $\beta$ service level (or fill rate) of a single-echelon supply chain system that tracks inventory using a reorder point policy.  The first model, Conventional, requires inventory to be tracked constantly, while the second model, Undershoot, represents a system where inventory is tracked periodically, which can potentially cause reordering delay thus exacerbating system backorders (called creating an undershoot).  The accuracy of these models  is tested using a discrete-event simulation built using SimPy.  We tested accuracy for different values of reorder quantity and demand coefficient of variation.  The models and simulation were compared for cases when system allows backorders (and late shipments) vs.\ the scenario when the unfulfilled demand is considered lost sales.

Based on our experiments, we make the following conclusions:
\bit
	\item In most cases tested in this work the closed form literature models under-predict the service level, and thus can result in a higher reorder point and safety stock
	\item Both Conventional and Undershoot models, for either backorder or lost sales cases, are valid predominantly when the reorder quantity is more than the reorder point (and should only be used when this is true)
	\item If reorder quantity is less than the reorder point and if the system follows a backorder strategy, the Conventional model can be used only when the coefficient of variation for demand is greater than 0.8
	\item In general, the systems following a lost sales approach should see higher $\beta$ service levels than the ones following a backorder strategy
\eit

\bibliographystyle{amsplain}

\providecommand{\bysame}{\leavevmode\hbox to3em{\hrulefill}\thinspace}
\providecommand{\MR}{\relax\ifhmode\unskip\space\fi MR }
\providecommand{\MRhref}[2]{%
  \href{http://www.ams.org/mathscinet-getitem?mr=#1}{#2}
}
\providecommand{\href}[2]{#2}

\end{document}